\documentclass[11pt]{article}
\usepackage{amsmath,amsthm,amsfonts,amssymb,enumerate,calc,graphicx,iwona}
\usepackage[colorlinks=true,citecolor=black,linkcolor=black,urlcolor=blue]{hyperref}
\usepackage[noabbrev,capitalise]{cleveref}
\usepackage[capitalise]{cleveref}
\usepackage[lmargin=30mm,rmargin=30mm,bmargin=30mm,tmargin=30mm]{geometry}
\usepackage[numbers,sort&compress]{natbib}
\renewcommand{\baselinestretch}{1.1}
\usepackage[hang]{footmisc} 

\renewcommand{\thefootnote}{\fnsymbol{footnote}}	
\allowdisplaybreaks
\newcommand\DateFootnote{
\begingroup
\renewcommand\thefootnote{}
\footnote{January 20, 2015, Revised \today}
\setcounter{footnote}{0}
\vspace*{-3ex}
\endgroup}

\makeatletter
\renewcommand\section{\@startsection {section}{1}{\z@}%
                                   {-3ex \@plus -1ex \@minus -.2ex}%
                                   {2ex \@plus.2ex}%
                                   {\normalfont\large\bfseries}}
\renewcommand\subsection{\@startsection{subsection}{2}{\z@}%
                                     {-2.5ex\@plus -1ex \@minus -.2ex}%
                                     {1.5ex \@plus .2ex}%
                                     {\normalfont\normalsize\bfseries}}
\renewcommand\subsubsection{\@startsection{subsubsection}{3}{\z@}%
                                     {-2ex\@plus -1ex \@minus -.2ex}%
                                     {1ex \@plus .2ex}%
                                     {\normalfont\normalsize\bfseries}}
 \renewcommand\paragraph{\@startsection{paragraph}{4}{\z@}%
                                    {1.5ex \@plus.5ex \@minus.2ex}%
                                    {-1em}%
                                    {\normalfont\normalsize\bfseries}}
\renewcommand\subparagraph{\@startsection{subparagraph}{5}{\parindent}%
                                       {1.5ex \@plus.5ex \@minus .2ex}%
                                       {-1em}%
                                      {\normalfont\normalsize\bfseries}}
\makeatother

\newcommand{\arXiv}[1]{arXiv:\,\href{http://arxiv.org/abs/#1}{#1}}
\newcommand{\msn}[1]{MR:\,\href{http://www.ams.org/mathscinet-getitem?mr=MR#1}{#1}}

\newcommand{\doi}[1]{doi:\,\href{http://dx.doi.org/#1}{#1}}

\newcommand{\Zbl}[1]{Zbl:\,\href{http://www.zentralblatt-math.org/zmath/en/search/?q=an:#1}{#1}}

\theoremstyle{plain}
\newtheorem{thm}{Theorem}
\newtheorem{lem}[thm]{Lemma}

\newtheorem{cor}[thm]{Corollary}

\theoremstyle{definition}

\newcommand{\FLOOR}[1]{\ensuremath{\protect\left\lfloor#1\right\rfloor}}

\newcommand{\ceil}[1]{\lceil{#1}\rceil}
\newcommand{\floor}[1]{\lfloor{#1}\rfloor}
\newcommand{\half}{\ensuremath{\protect\tfrac{1}{2}}}
\renewcommand{\geq}{\geqslant}
\renewcommand{\leq}{\leqslant}

\newcommand{\BY}{Bourgain and Yehudayoff~\citep{Bourgain09,BY13}}

\begin{document}


{\Large\bfseries\boldmath\scshape Layouts of Expander Graphs}

\DateFootnote

Vida Dujmovi{\'c}\,\footnotemark[2]
\quad 
Anastasios Sidiropoulos\,\footnotemark[3]
\quad  
David~R.~Wood\,\footnotemark[4]

\footnotetext[2]{School of Computer Science and Electrical Engineering, University of Ottawa, Ottawa, Canada (\texttt{vida.dujmovic@uottawa.ca}).   Supported by NSERC.} 

\footnotetext[3]{Department of Computer Science and Engineering, and Department of Mathematics, Ohio State University (\texttt{sidiropo@gmail.com}).} 

\footnotetext[4]{School of Mathematical Sciences, Monash University, Melbourne, Australia (\texttt{david.wood@monash.edu}). Research supported by  the Australian Research Council.}

\emph{Abstract.} Bourgain and Yehudayoff recently constructed $O(1)$-monotone bipartite expanders. By combining this result with a generalisation of the unraveling method of Kannan, we construct 3-monotone bipartite expanders, which is best possible. We then show that the same graphs admit 3-page book embeddings, 2-queue layouts, 4-track layouts, and have simple thickness 2. All these results are best possible. 


\section{Introduction}
\label{Intro}

\renewcommand{\thefootnote}{\arabic{footnote}}

Expanders are classes of highly connected graphs that are of fundamental importance in graph theory, with numerous applications, especially in theoretical computer science \citep{HLW-BAMS06}. While the literature contains various definitions of expanders, this paper focuses on bipartite expanders. For $\epsilon\in(0,1]$, a bipartite graph $G$ with bipartition $V(G)=A\cup B$ is a \emph{bipartite $\epsilon$-expander} if $|A|=|B|$ and $|N(S)|\geq(1+\epsilon)|S|$ for every subset $S\subset A$ with $|S|\leq\frac{|A|}{2}$. Here $N(S)$ is the set of vertices adjacent to some vertex in $S$. An infinite family of bipartite $\epsilon$-expanders, for some fixed $\epsilon>0$, is called an \emph{infinite family of bipartite expanders}. 

There has been much research on constructing and proving the existence of expanders with various desirable properties. The first example is that there is an infinite family of expanders with bounded degree, in fact, degree at most 3 (see \citep{ASS08,RVV-AM02,HLW-BAMS06} for example). 

\subsection{Monotone Layouts}

\BY\ recently gave an explicit construction of an infinite family of bipartite expanders with an interesting additional property. Say $G$ is a bipartite graph with ordered colour classes $(v_1,\dots,v_n)$ and $(w_1,\dots,w_m)$. Two edges $v_iw_j$ and $v_kw_\ell$ \emph{cross} if $i<k$ and $\ell<j$. A matching $M$ in $G$ is \emph{monotone} if no two edges in $M$ cross. A bipartite graph with ordered coloured classes is \emph{$d$-monotone} if it is the  union of $d$ monotone matchings. Note that every $d$-monotone bipartite graph has maximum degree at most $d$. Motivated by connections to dimension expanders, \citet{DS11} constructed an infinite family of $O(\log n)$-monotone bipartite expanders\footnote{While monotone expanders are not explicitly mentioned by  \citet{DS11}, the connection is made explicit by \citet{DW-ToC10}.}. \citet{DW-ToC10}  constructed an infinite family of $O(\log^c n)$-monotone bipartite expanders, for any constant $c>0$. \BY\ proved the following breakthrough\footnote{An outline of the proof was given in the original paper by \citet{Bourgain09}. A full proof was given by \citet{BY13}. See the paper by \citet{DW-ToC10} for more discussion.}:

\begin{thm}[\BY]
\label{Bourgain}
There is an infinite family of $d$-monotone bipartite expanders, for some constant $d$. 
\end{thm}

Note that the proof of \cref{Bourgain} is constructive, and indeed no probabilistic proof is known. This is unusual, since probabilistic proofs for the existence of expanders are typically easier to obtain than explicit constructions. 

The first contribution of this paper is to  show how any $O(1)$-monotone bipartite expander can be manipulated to produce a $3$-monotone bipartite expander.

\begin{thm}
\label{3Monotone}
There is an infinite family of $3$-monotone bipartite expanders.
\end{thm}

\subsection{Book Embeddings}

\cref{3Monotone} has applications to \emph{book embeddings}. A \emph{$k$-page book embedding} of a graph $G$ consists of a linear order $(u_1,\dots,u_n)$ of $V(G)$ and a partition $E_1,\dots,E_k$ of $E(G)$, such that edges in each set $E_i$ do not cross with respect to $(u_1,\dots,u_n)$. That is, for all $i\in[1,k]$, there are no edges $u_au_b$ and $u_cu_d$ in $E_i$ with $a<c<b<d$. One may think of the vertices as being ordered along the spine of a book, with each edge drawn on one of $k$ pages, such that no two edges on the same page cross. A graph with a $k$-page book embedding is called a \emph{$k$-page graph}. The \emph{page-number} of a graph $G$ is the minimum integer $k$ such that there is a $k$-page book embedding of $G$. Note that page-number is also called \emph{book thickness} or \emph{stack-number} or \emph{fixed outer-thickness}; see reference \citep{DujWoo-DMTCS04} for more on book embeddings. A $k$-page book embedding is \emph{$k$-pushdown} if, in addition, each set $E_i$ is a matching \citep{GKS89}. 

\begin{figure}[!htb]
\begin{center}
\includegraphics{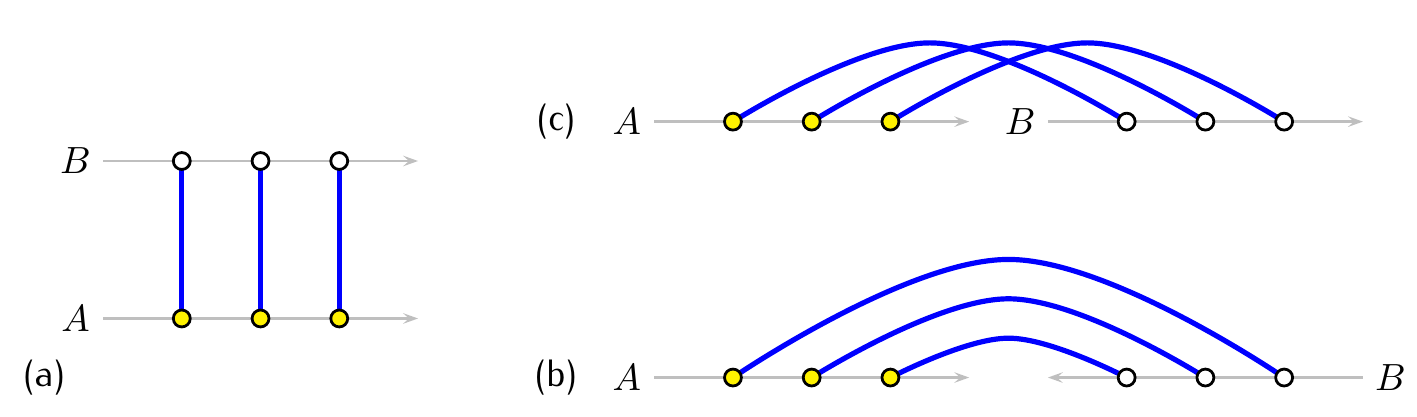}
\caption{Converting (a) a monotone matching to (b) a book embedding and (c) a queue layout~\citep{DW-ToC10,DPW04,Pemmaraju-PhD}.}
\label{TrackStackQueue}
\end{center}
\end{figure}

A $d$-monotone bipartite graph has a $d$-pushdown book embedding, and thus has page-number at most $d$, since using the above notation, edges in a monotone matching do not cross in the vertex ordering $(v_1,\dots,v_n,w_m,\dots,w_1)$, as illustrated in \cref{TrackStackQueue}(b). This observation has been made several times in the literature \citep{DW-ToC10,DPW04,Pemmaraju-PhD}. In the language of \citet{Pemmaraju-PhD}, this book embedding is `separated'.

Thus the construction of \BY\ gives an infinite family of $d$-pushdown bipartite expanders with maximum degree $d$, for some constant $d$. This result solves an old open problem of \citet{GKS89,GKS-JCSS89} that arose in the  modelling of multi-tape Turing machines. In particular, \citet{GKS-JCSS89} showed that there are $O(1)$-pushdown expanders if and only if it is not possible for a 1-tape nondeterministic Turing machine to simulate a 2-tape machine in subquadratic time.

\cref{3Monotone} and the above observation implies:

\begin{thm}
\label{3Page}
There is an infinite family of $3$-pushdown bipartite expanders. 
\end{thm}

\subsection{Queue Layouts}

Queue layouts are dual to book embeddings. (In this setting, book embeddings are often called stack layouts.)\ A \emph{$k$-queue layout} of a graph $G$ consists of a linear order $(u_1,\dots,u_n)$ of $V(G)$ and a partition $E_1,\dots,E_k$ of $E(G)$, such that edges in each set $E_i$ do not \emph{nest} with respect to $(u_1,\dots,u_n)$. That is, for all $i\in[1,k]$, there are no edges $u_au_b$ and $u_cu_d$ in $E_i$ with $a<c<d<b$. A graph with a $k$-queue layout is called a \emph{$k$-queue graph}. The \emph{queue-number} of a graph $G$ is the minimum integer $k$ such that there is a $k$-queue layout of $G$. See \citep{DMW05,HLR92,HR92,DujWoo-DMTCS04,DPW04,DFP13,DMW13,DujWoo-DMTCS05} and the references therein for results on queue layouts. 

A $d$-monotone bipartite graph has queue-number at most $d$, since using the above notation, edges in a monotone matching do not cross in the vertex ordering $(v_1,\dots,v_n,w_1,\dots,w_m)$, as illustrated in \cref{TrackStackQueue}(c).  Thus the construction of \BY\ provides an infinite family of bipartite expanders with bounded queue-number, as observed by \citet{DMS14}. And \cref{3Monotone} gives an infinite family of 3-queue bipartite expanders. We improve this result as follows.

\begin{thm}
\label{2Queue}
There is an infinite family of $2$-queue bipartite expanders with maximum degree $3$.
\end{thm}

\subsection{Track Layouts}

Finally, consider track layouts of graphs. In a graph $G$, a \emph{track} is an independent set, equipped with a total ordering denoted by $\preceq$.  A \emph{$k$-track layout} of a graph $G$ consists of a partition $(V_1,\dots,V_k)$ of $V(G)$ into tracks, such that between each pair of tracks, no two edges cross. That is, there are no edges $vw$ and $xy$ in $G$ with $v\prec x$ in some track $V_i$, and $y\prec w$ in some track $V_j$. The \emph{track-number} is the minimum integer $k$ for which there is a $k$-track layout of $G$. See \citep{DMW05,DPW04,DMW13,DLMW-DM09,DujWoo-DMTCS05} and the references therein for results on track layouts.  We prove the following:


\begin{thm}
\label{4Track}
There is an infinite family of $4$-track bipartite expanders with maximum degree $3$.
\end{thm}

\subsection{Discussion}

Some notes on the above theorems are in order. First note that the proofs of Theorems~\ref{3Monotone}--\ref{4Track} are unified. Indeed, each of these theorems refer to the same family of graphs.

\paragraph{Tightness:}  
Each of Theorems~\ref{3Monotone}--\ref{4Track} is best possible since 2-page graphs (and thus 2-monotone graphs) are planar \cite{BK79}, 1-queue graphs are planar \cite{HR92}, and 3-track graphs are planar \cite{DujWoo-DMTCS05}, but planar graphs have $O(\sqrt{n})$ separators \citep{LT79}, and are thus far from being expanders. It is interesting that graphs that are `close' to being planar can be expanders. 


\paragraph{Expansion and Separators:} \citet{Sparsity} introduced the following definition. A class $\mathcal{G}$ of graphs has \emph{bounded expansion} if there is a function $f$ such that for every integer $r\geq 0$ and every graph $G\in\mathcal{G}$, any graph obtained from $G$ by contracting disjoint balls of radius $r$ has average degree at most $f(r)$. The least such function $f$ is called the \emph{expansion function} for $\mathcal{G}$. For example, minor-closed classes have constant expansion functions (independent of $r$). \citet{NOW} proved that graph classes with bounded page-number or bounded queue-number have bounded expansion (also see \citep[Chapter~14]{Sparsity}).  Thus, Theorems~\ref{3Page} and \ref{2Queue}  provide natural families of graphs that have bounded expansion yet contain an infinite family of expanders. The upper bound (proved in \citep{NOW}) on the expansion function for graphs of bounded page-number or bounded queue-number is exponential. \citet{Sparsity} state as an open problem whether this exponential bound is necessary. Since graph classes with sub-exponential expansion functions have $o(n)$ separators \citep[Theorem~8.3]{NesOdM-GradII} (also see \citep{Dvorak14}), and expanders do not have $o(n)$ separators (see \cref{Separators}), Theorems~\ref{3Page} and \ref{2Queue} imply that  indeed exponential expansion is necessary for 3-page and 2-queue graphs. Since queue-number is tied to track-number \citep{DPW04}, these same conclusions hold for track-number. 

\paragraph{Subdivisions:} Theorems  slightly weaker than Theorems~\ref{3Page}--\ref{4Track} can be proved using subdivisions. It can be proved that if $G$ is a bipartite $\epsilon$-expander with bounded degree, then the graph obtained from $G$ by subdividing each edge twice is a bipartite $\epsilon'$-expander (see  \cref{SubdivExpander}). \citet{DujWoo-DMTCS05} proved that every $k$-page graph has a $3$-page subdivision with $2\ceil{\log_2 k} -2$ division vertices per edge. Applying this result to the construction of \BY, we obtain an infinite family of 3-page bipartite expanders with bounded degree. Note that the degree bound here is the original degree bound from the construction \BY, which is much more than 3 (the degree bound in \cref{3Page}). In particular, 3-monotone expanders cannot be constructed using subdivisions. 

One can also construct 2-queue expanders and 4-track expanders using subdivisions. \citet{DujWoo-DMTCS05} proved that every $k$-queue graph has a $2$-queue subdivision with $2\ceil{\log_2 k} + 1$ division vertices per edge, and has a $4$-track subdivision with $2\ceil{\log_2 k} + 1$ division vertices per edge. To apply these results, one must modify the relevant constructions so that each edge is subdivided an even number of times (details omitted). Again the obtained degree bound is weaker than in Theorems~\ref{2Queue} and \ref{4Track}. 

\paragraph{Thickness:}  The \emph{thickness} of a graph $G$ is the minimum integer $k$ such that $G=G_1\cup\dots\cup G_k$ for some planar subgraphs $G_1,\dots,G_k$. See \citep{MOS98} for a survey on thickness. A natural question arises: what is the minimum integer $k$ for which there is an infinite family of bipartite expanders with thickness $k$? It is easily seen that there are bipartite expanders with thickness 2: Let $G'$ be the graph obtained from an $\epsilon$-bipartite expander $G$ with bounded degree by subdividing each edge twice. Then $G'$ is an $\epsilon'$-expander (see \cref{SubdivExpander}). The edges of $G'$ incident to the original vertices form a star forest $G_1$, and the remaining edges form a matching $G_2$, both of which are planar. Hence $G'$ has thickness 2. Of course, thickness 2 is best possible for an expander since every graph with thickness 1 is planar. 

Every graph with thickness $k$ can be drawn in the plane with no crossings between edges in each of the $k$ given planar subgraphs (since a planar graph can be drawn without crossings with its vertices  at prespecified positions). However edges from different planar subgraphs might cross multiple times. This motivates the following definition. A drawing of a graph is \emph{simple} if no two edges cross more than once. The \emph{simple thickness} of a graph $G$ is the minimum integer $k$ such that there is a simple drawing of $G$ and a partition of $E(G)$ into $k$ non-crossing subgraphs. 

\begin{figure}[!ht]
\begin{center}
\includegraphics{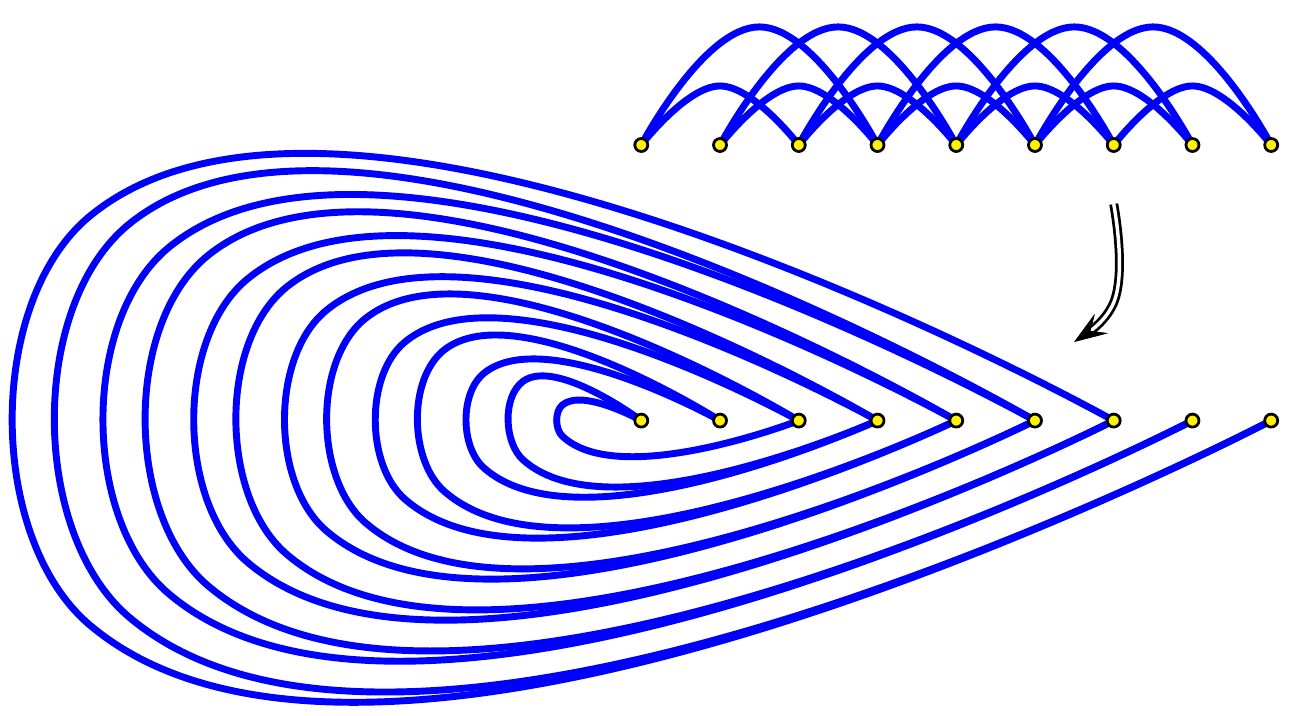}
\caption{Drawing a 1-queue graph without crossings.}
\label{PlanarQueue}
\end{center}
\end{figure}

We now show how to obtain an infinite family of bipartite graphs with simple thickness 2. Every 1-queue graph is planar  \cite{HR92}. To see this, say $v_1,\dots,v_n$ is the vertex ordering in a 1-queue graph. Position $v_i$ at $(i,0)$ in the plane. Draw each edge $v_iv_j$ with $i<j$, as a curve from $(i,0)$ starting above the X-axis, through $(-i-j,0)$, and then under the X-axis to $(j,0)$, as illustrated in \cref{PlanarQueue}. Since no two edges are nested in the initial ordering, no two edges cross in this drawing. Now, given a 2-queue layout, applying the same construction for each queue gives a simple drawing, in which edges from the first queue do not cross, edges from the second queue do not cross, and each edge from the first queue crosses each edge from the second queue at most once (if the curves are drawn carefully). This shows that every 2-queue graph has a simple drawing with thickness 2. By \cref{2Queue} there is an infinite family of bipartite expanders with simple thickness 2.  Furthermore, one may subdivide each edge twice in the above construction, and then draw each edge straight to obtain an infinite family of bipartite expanders with geometric thickness 2 (see \citep{BMW-EJC06,DEH-JGAA00,DujWoo-DCG07}). 
 

\section{Two-Sided Bipartite Expanders}
\label{TwoSided}

Throughout this paper, it is convenient to employ the following definition. A bipartite graph $G$ with bipartition $A,B$ is a \emph{two-sided bipartite $\epsilon$-expander} if $|A|=|B|$, and for all $S\subset A$ with $|S|\leq\frac{|A|}{2}$ we have $|N(S)|\geq(1+\epsilon)|S|$, and  for all $T\subset B$ with $|T|\leq\frac{|B|}{2}$ we have $|N(T)|\geq(1+\epsilon)|T|$. This is a strengthening of the notion of a (one-sided) bipartite $\epsilon$-expander. The next lemma says that a (one-sided) $k$-monotone bipartite expander can be easily modified to produce a two-sided $2k$-monotone bipartite expander. 

\begin{lem}
\label{TwoSIded}
If $G$ is a $k$-monotone bipartite $\epsilon$-expander with ordered bipartition $A=(v_1,\dots,v_n)$ and $B=(w_1,\dots,w_m)$, then the graph $G'$ with vertex set $V(G'):=V(G)$ and edge set $E(G'):=\{(v_i,w_j),(v_j,w_i):(v_i,w_j)\in E(G)\}$ is a two-sided $2k$-monotone  bipartite $\epsilon$-expander. 
\end{lem}

\begin{proof}
Observe that $(v_i,w_j)$ crosses $(v_a,w_b)$ if and only if $(v_j,w_i)$ crosses $(v_b,w_a)$. Thus if $M$ is a monotone matching, then $\{(v_j,w_i):(v_i,w_j)\in M\}$ is also a monotone matching. Hence, $E(G')$ can be partitioned into $2k$ monotone matchings. Since $G$ is a spanning subgraph of $G'$,  we have that $G'$ is an $\epsilon$-expander. Given $T\subseteq B$ with $|T|\leq\frac{|B|}{2}$, define $S:=\{v_i\in A:w_i\in B\}$. Then $|N_{G'}(T)|\geq|N_G(S)|\geq(1+\epsilon)|S|=(1+\epsilon)|T|$. Thus $G'$ is a two-sided bipartite $\epsilon$-expander. 
\end{proof}

The construction of \BY\ and \cref{TwoSIded} together imply:

\begin{cor}
\label{TwoSidedBourgain}
There is an infinite family of two-sided $d$-monotone bipartite expanders, for some constant $d$.
\end{cor}

\section{Unraveling}
\label{Unraveling}

\begin{figure}[!ht]
\begin{center}
\includegraphics{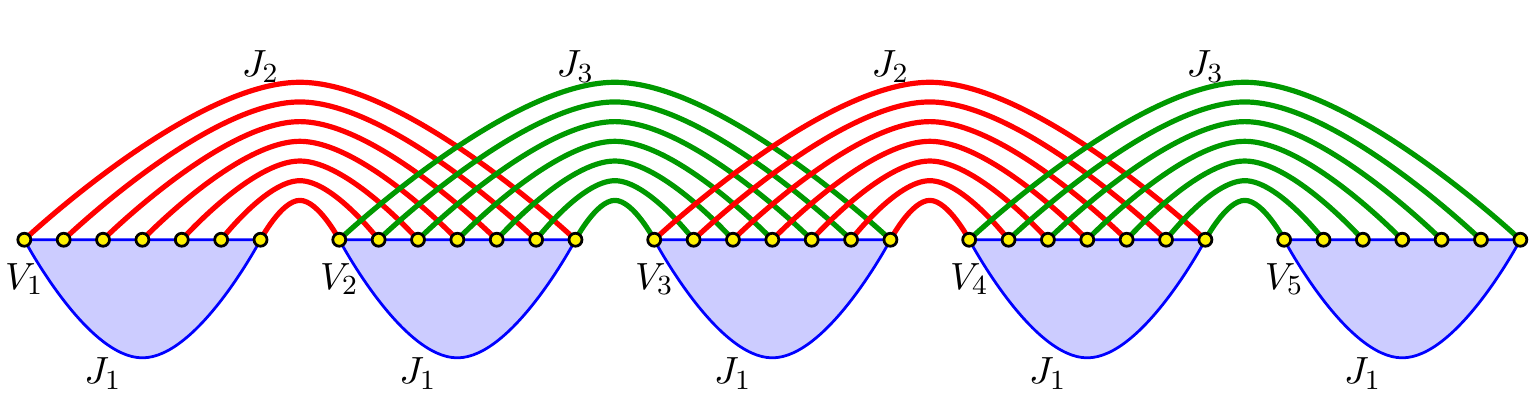}
\caption{3-page book embedding of the unraveling, due to \citet{Kannan85}.}
\label{Kannan}
\end{center}
\end{figure}

The following construction of \citet{Kannan85} is the starting point for our work. Let $G$ be a graph, whose edges are $k$-coloured (not necessarily properly). Let $E_1,\dots,E_k$ be the corresponding partition of $E(G)$. Let $G'$ be the graph with vertex set $$V(G'):=V(G)\times[1,k]=\{v_i:v\in V(G),i\in[1,k]\},$$ where $v_iw_i\in E(G')$ for each edge $vw\in E_i$ and $i\in[1,k]$, and $v_iv_{i+1}\in E(G')$ for each vertex $v\in V(G)$ and $i\in[1,k-1]$. \citet{Kannan85} called $G'$ the \emph{unraveling} of $G$, which he defined  in the case that the edge colouring comes from a $k$-page book embedding, and proved that $G'$ has a 3-page book embedding. To see this, for $i\in[1,k]$, let $V_i:=\{v_i:v\in V(G)\}$ ordered by the given ordering of $V(G)$. Define 
\begin{align*}
J_1&:=\{v_iw_i:vw\in E_i,i\in[1,k]\}\\
J_2&:=\{v_iv_{i+1}:v\in V(G),i\in[1,k], i\text{ odd} \} \text{ and }\\  
J_3&:=\{v_iv_{i+1}:v\in V(G),i\in[1,k], i\text{ even} \}.
\end{align*} 
Then $J_1,J_2,J_3$ is a partition of $E(G')$, and for $i\in[3]$, no two edges in $J_i$ cross with respect to the vertex ordering $V_1,V_2,\dots,V_k$, as illustrated in \cref{Kannan}. Thus, this is a 3-page book embedding of $G'$.

This observation is extended as follows. 

\begin{lem}
\label{3MonotoneUnraveling}
If a bipartite graph $G$ is $k$-monotone, then the unraveling $G'$ is 3-monotone.
\end{lem}

\begin{figure}[!ht]
\begin{center}
\includegraphics{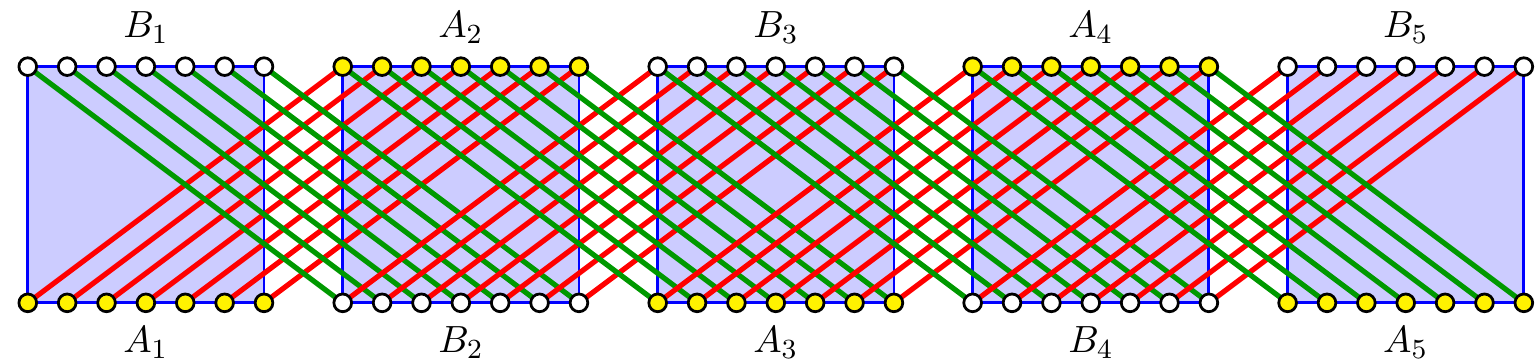}
\caption{3-monotone layout of the unraveling.}
\label{3MonotoneUnravelling}
\end{center}
\end{figure}

\begin{proof}
Say $A,B$ is the given bipartition of $G$. Let $A_i:=V_i\cap A$ and $B_i:=V_i\cap B$, where $V_i$ is defined above. Say $A_i$ and $B_i$ inherit the given orderings of $A$ and $B$ respectively. Then $G'$ is bipartite with ordered bipartition given by $A_1,B_2,A_3,B_4,A_5,B_6,\dots$ and $B_1,A_2,B_3,A_4,B_5,A_6\dots$. Observe that for $i\in[3]$, no two edges in $J_i$ cross with respect to these orderings, as illustrated in \cref{3MonotoneUnravelling}. Thus $G'$ is 3-monotone.
\end{proof}

The unraveling $G'$ has interesting expansion properties. In particular, \citet{Kannan85} proved that if $G'$ has a small separator, then so does $G$. Thus, if $G$ is an expander, then $G$ and $G'$ have no small separator. Various results in the literature say that if every separator of an $n$-vertex graph $G$ has size at least $\epsilon n$, then $G$ contains an expander as a subgraph (for various notions of non-bipartite expansion). However, the unraveling $G'$ might not be a bipartite expander. For example, $G'$ might have a vertex of degree 1. This happens for a vertex $v_1$ where $v$ is incident to no edge coloured $1$, or a vertex $v_k$ where $v$ is incident to no edge coloured $k$. The natural solution for this problem is to add the edge $v_1v_k$ for each vertex $v$ of $G$. Now each vertex $v$ corresponds to the cycle $C_v=(v_1,v_2,\dots,v_k)$. However, the obtained graph is still not an expander: if $S$ consists of every second vertex in some $C_v$, then it is possible for $N(S)$ to consist only of the other vertices in $C_v$, in which case $|N(S)|=|S|$, and the graph is not an expander. Moreover, it is far from clear how to construct a 3-monotone layout of this graph. (For even $k$, the layout in the proof of \cref{3MonotoneUnraveling} is 5-monotone.)\ 

\section{Generalised Unraveling}

The obstacles discussed at the end of the previous section are overcome in the following  lemma. This result is reminiscent of the replacement product; see \citep{ASS08,RVV-AM02,HLW-BAMS06,DW-ToC10}.

\begin{lem}
\label{Gen}
Let $G$ be a two-sided bipartite $\epsilon$-expander with bipartition $A,B$ and maximum degree $\Delta$. 
Let $n:=|A|=|B|$. Assume $n\geq3$. Let $k\geq 2$ be an integer. For each vertex $v$ of $G$, let $k_v$ be an integer with $k\leq k_v\leq (1+\frac{\epsilon}{4})k$.  Let $G'$ be a bipartite graph with bipartition $X,Y$ such that:
\begin{itemize}
\item $G'$ contains disjoint cycles $\{C_v:v\in V(G)\}$, 
\item $|C_v|=2k_v$ for each vertex $v\in V(G)$, 
\item $V(G')=\cup\{V(C_v):v\in V(G)\}$,  and
\item for each edge $vw$ of $G$ there are edges $xy$ and $pq$ of $G'$ such that $x\in C_v\cap X$ and $y\in C_w\cap Y$ and $p\in C_v\cap Y$ and $q\in C_w\cap X$. 
\end{itemize}
Then $G'$ is a two-sided bipartite $\epsilon'$-expander, for some $\epsilon'$ depending only on $\epsilon$, $k$ and $\Delta$. 
\end{lem}

\begin{proof}
For each vertex $v$ of $G$, we have $|C_v\cap X|=|C_v\cap Y|=k_v$. Thus 
$$|X|=|Y|=\sum_{v\in V(G)}k_v\leq 2(1+\tfrac{\epsilon}{4})kn.$$ 
Let $S\subseteq X$ with $|S|\leq\frac{|X|}{2}$, which is at most $(1+\frac{\epsilon}{4})kn$. 
By the symmetry between $X$ and $Y$, it suffices to prove that $|N_{G'}(S)|\geq (1+\epsilon')|S|$. 

For each vertex $v$ of $G$, observe that $|C_v\cap S|\leq|C_v\cap X|=k_v$.
Say $v$ is \emph{heavy} if $|C_v\cap S|=k_v$.
Say $v$ is \emph{light} if $1\leq |C_v\cap S|\leq k_v-1$.
Say $v$ is \emph{unused} if $C_v\cap S=\emptyset$.
Each vertex of $G$ is either heavy, light or unused. 

Say a heavy vertex $v$ of $G$ is \emph{fat} if every neighbour of $v$ is also heavy.
Let $F$ be the set of fat vertices in $G$. 
Let $H$ be the set of non-fat heavy vertices in $G$. 
Let $L$ be the set of light vertices in $G$. 
Let $U$ be the set of unused vertices in $G$. 
Thus $F,H,L,U$ is a partition of $V(G)$. 
Let $f:=|F|$ and $h:=|H|$ and $\ell:=|L|$. 
Let $f_A:=|F\cap A|$ and $f_B:=|F\cap B|$ and $h_A:=|H\cap A|$ and $h_B:=|H\cap B|$. 

Since the vertices in $H$ are not fat, every vertex in $H$ has a neighbour in $L\cup U$. Let $H'$ be the set of vertices in $H$ adjacent to no vertex in $U$ (and thus with a neighbour in $L$). Let $H''$ be the set of vertices in $H$ adjacent to some vertex in $U$. Define $h':=|H'|$ and $h'':=|H''|$. 

For each vertex $v$ of $G$, let $c(v)$ be the number of vertices in $C_v$ adjacent to some vertex in $C_v\cap S$. Since $C_v\cap S$ is an independent set in $C_v$, by \cref{IndSetCycle} below, if $v$ is heavy, then every second vertex of $C_v$ is in $S$ and $c(v)=k_v=|C_v\cap S|$, and if $v$ is light then $c(v)\geq |C_v\cap S|+1$. Thus
\begin{equation}
\label{FirstCount}
|N_{G'}(S)| \;\geq\; 
\sum_{v\in F\cup H\cup L} \!\!\!\! c(v) \;\geq\; 
\ell + \!\!\!\!\sum_{v\in F\cup H\cup L} \!\!\!\! |C_v\cap S| \;=\; 
\ell + |S|. 
\end{equation}
Moreover, each vertex $v$ in $H''$ is adjacent in $G$ to some vertex $w$ in $U$. By assumption, there is an edge $xy$ of $G'$ such that $x\in C_v\cap X$ and $y\in C_w\cap Y$.  Since $v$ is heavy, $x$ is in $S$ and $y$ is in $N_{G'}(S)$. And since $w$ is unused, $y$ is adjacent to no vertex in $C_w\cap S$. Thus $y$ is not counted in the lower bound on $N_{G'}(S)$ in \eqref{FirstCount}. Each such vertex $y$ is adjacent to at most $\Delta$ vertices in $H''$. Hence
\begin{equation}
\label{SecondCount}
|N_{G'}(S)| \geq |S| + \ell + \frac{h''}{\Delta}.
\end{equation}
Our goal now is to prove that $\ell+\frac{h''}{\Delta} \geq\epsilon'|S|$, where $\epsilon':=(k+k\Delta(\frac{1+\epsilon}{\epsilon}))^{-1}$. 
Since
$1-\epsilon'(k+k\Delta(1+\tfrac{1}{\epsilon})= 0$ and
$\tfrac{1}{\Delta}-\epsilon' k(1+\tfrac{1}{\epsilon})\geq 0$,
$$\big(1-\epsilon'(k+k\Delta(1+\tfrac{1}{\epsilon})\big)\ell + 
\big(\tfrac{1}{\Delta}-\epsilon' k(1+\tfrac{1}{\epsilon})\big)h'' \geq 0.$$
That is, 
$$\ell + \frac{h''}{\Delta} 
\geq 
\epsilon'k\ell + \epsilon' k(1+\tfrac{1}{\epsilon})(\Delta\ell+h'') .$$
Every vertex in $H'$ has a neighbour in $L$, each of which has degree at most $\Delta$.  Thus $h'\leq \Delta \ell$, implying
$$\ell +  \frac{h''}{\Delta}
\geq 
\epsilon' k\ell + 
\epsilon' k(1+\tfrac{1}{\epsilon})(h'+h'') 
=
\epsilon' k(\ell + h+\tfrac{h}{\epsilon})
.$$
Suppose, on the contrary, that $f_A\geq\floor{\frac{n}{2}}+1$. Let $Q$ be a subset of $F\cap A$ of size $\floor{\frac{n}{2}}$. Since $G$ is a two-sided $\epsilon$-expander, and since every neighbour of each vertex in $F\cap A$ is in $(F\cup H)\cap B$, we have
$f_B+h_B\geq |N_G(Q)|\geq(1+\epsilon)\floor{\frac{n}{2}}$. Thus 
$$f_A+f_B+h_B\geq \floor{\tfrac{n}{2}}+1+(1+\epsilon)\floor{\tfrac{n}{2}}
\geq 
(2+\epsilon)\floor{\tfrac{n}{2}}+1
\geq 
n +\tfrac{\epsilon}{2}(n-1)
.$$ 
However, 
$(1+\frac{\epsilon}{4})kn\geq |S|\geq k(f_A+f_B+h_B)$, implying $f_A+f_B+h_B\leq (1+\frac{\epsilon}{4})n$, which is a contradiction (since $n\geq 3$). 

Now assume that $f_A\leq \frac{n}{2}$. Since $G$ is a two-sided $\epsilon$-expander, and since every neighbour of a vertex in $F\cap A$ is in $(F\cup H)\cap B$,  we have
$f_B+h_B\geq |N_G(F\cap A)|\geq(1+\epsilon)f_A$. 
By symmetry, $f_A+h_A\geq (1+\epsilon)f_B$. 
Thus $f+h\geq(1+\epsilon)f$ and $h\geq\epsilon f$. Hence
$$\ell +  \frac{h''}{\Delta}
\geq 
\epsilon' k(\ell+h+f)
.$$
Since $|S|\leq k(f+h+\ell)$, 
$$
\ell  +  \frac{h''}{\Delta}
\geq 
\epsilon' |S|,$$
and
$N_{G'}(S)\geq(1+\epsilon')|S|$, as desired. 
\end{proof}

\begin{lem}
Let $I$ be an independent set in a a cycle graph $C$. Then $|N_C(I)|\geq|I|$ with equality only if $I=\emptyset$ or $|C|=2|I|$. 
\label{IndSetCycle}
\end{lem}

\begin{proof}
For each vertex $x$ in $N_C(I)$, if $x$ is adjacent to exactly one vertex $v$ in $I$, then send the charge of 1 from $x$ to $v$, and if $x$ is adjacent to exactly two vertices $v$ and $w$ in $I$, then send a charge of $\tfrac12$ from $x$ to each of $v$ and $w$. Each vertex in $I$ receives a charge of at least $\frac12$ from each of its neighbours in $C$. Thus the total charge, $|N_C(I)|$, is at least $I$, as claimed. If the total charge equals $|I|$, then each vertex $v$ in $I$ receives a charge of exactly 1, which implies that both neighbours of $v$ sent a charge of $\frac12$ to $v$. Thus both neighbours of $v$ are adjacent to two vertices in $I$. It follows that $I$ consists of every second vertex in $C$, and $|C|=2|I|$. 
\end{proof}

\section{The Wall}
\label{Walls}

\begin{figure}[h]
\begin{center}
\includegraphics{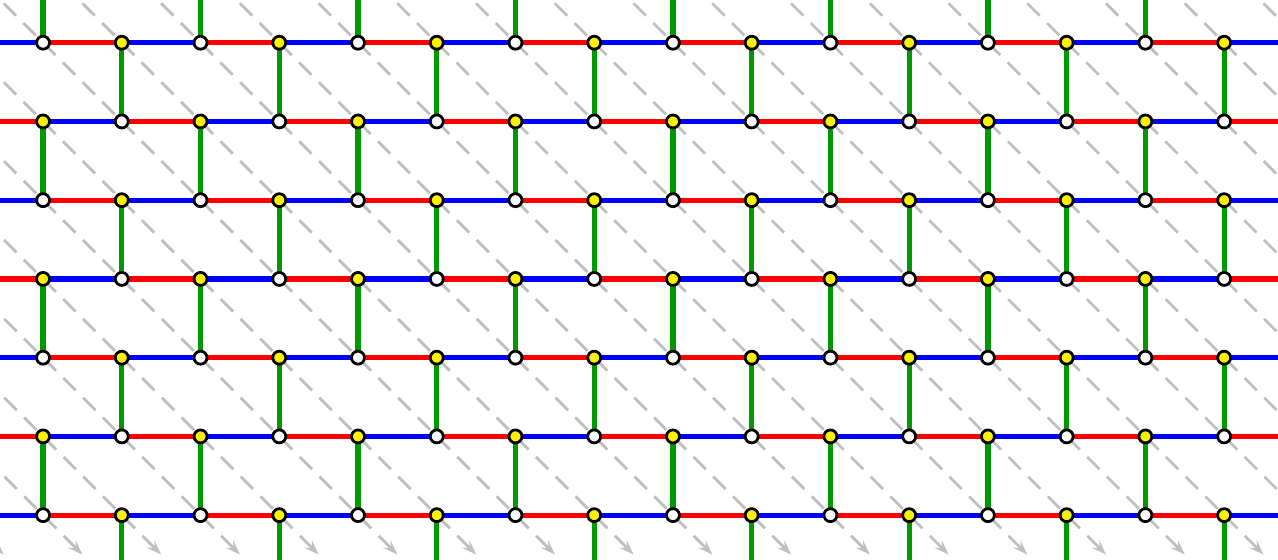}
\caption{The wall with ordered colour classes.\label{wall}}
\end{center}
\end{figure}

The following example is a key to our main proofs, and is of independent interest. The \emph{wall} is the infinite graph $W$ with vertex set $\mathbb{Z}^2$ and edge set 
$$\big\{\{(x,y)(x+1,y)\}:x,y\in\mathbb{Z}\big\}\cup\big\{\{(x,y)(x,y+1)\}:x,y\in\mathbb{Z}^+,\,x+y\text{ even}\big\}.$$
As illustrated in \cref{wall}, the wall is 3-regular and planar.

The next two results depend on the following vertex ordering of $W$. For vertices $(x,y)$ and $(x',y')$ of $W$, define 
$(x,y)\preceq (x',y')$ if $x+y<x'+y'$, or $x+y=x'+y'$ and $x\leq x'$.


\begin{lem}
\label{Wall}
The wall is 3-monotone bipartite.
\end{lem}

\begin{proof}
Let $A:=\{(x,y)\in\mathbb{Z}^2:\,x+y\text{ even}\}$ and  $B:=\{(x,y)\in\mathbb{Z}^2:\,x+y\text{ odd}\}$. 
Observe that $A,B$ is a bipartition of $W$. 
Consider $A$ and $B$ to be ordered by $\preceq$. 
Colour the edges of $W$ as follows. For each vertex $(x,y)$ where $x+y$ is even, 
colour $(x,y)(x+1,y)$ red, 
colour $(x,y)(x-1,y)$ blue, and
colour $(x,y)(x,y+1)$ green, as illustrated in \cref{wall}. 
Each edge of $W$ is thus coloured. 
If $(x,y)\prec (x',y')$ in $A$, then $(x+1,y)\prec(x'+1,y)$ in $B$. 
Thus the red edges form a monotone matching. 
Similarly, the green edges form a monotone matching, 
and the blue edges form a monotone matching. 
Thus $W$ is 3-monotone.
\end{proof}


\begin{lem}
\label{WallQueue}
The wall has a 2-queue layout, such that for all edges $pq$ and $pr$ with $p\prec q\prec r$ or $r\prec q\prec p$, the edges $pq$ and $pr$ are in distinct queues (called a `strict' 2-queue layout in \citep{Wood-Queue-DMTCS05}).
\end{lem}

\begin{proof}
We first prove that no two edges of $W$ are nested with respect to $\preceq$. 
Suppose that some edge $(x_2,y_2)(x_3,y_3)$ is nested inside another edge $(x_1,y_1)(x_4,y_4)$, 
where $(x_1,y_1)\prec(x_2,y_2)\prec(x_3,y_3)\prec(x_4,y_4)$. 
By the definition of $\preceq$, we have
$x_1+y_1\leq x_2+y_2\leq x_3+y_3\leq x_4+y_4$. Since
$(x_2,y_2)(x_3,y_3)$ and $(x_1,y_1)(x_4,y_4)$ are edges, 
$x_4+y_4=x_1+y_1+1$ and $x_3+y_3=x_2+y_2+1$. 
Hence $x_1+y_1= x_2+y_2$ and $x_3+y_3= x_4+y_4$. 
By the definition of $\preceq$, we have 
$x_1<x_2$ and $y_2<y_1$, and $x_3<x_4$ and $y_4<y_3$. 
Since $(x_1,y_1)(x_4,y_4)$ is an edge with  $(x_1,y_1)\prec(x_4,y_4)$, either $x_4=x_1+1$ or $y_4=y_1+1$. 
First suppose that $x_4=x_1+1$. Then $x_1<x_2\leq x_3$, implying $x_4=x_1+1\leq x_3$, which is a contradiction. 
Now assume that $y_4=y_1+1$. Then $y_1+1=y_4<y_3$. 
Since $(x_2,y_2)(x_3,y_3)$ is an edge, $y_3 \leq y_2+1$, implying $y_1<y_2$, which is a contradiction. 
Hence no two edges are nested. 

For each vertex $(x,y)$ where $x+y$ is even, assign the edges $(x,y)(x+1,y)$ and $(x,y)(x-1,y)$ to the first queue, and assign the edge $(x,y)(x,y+1)$ to the second queue. If $x+y$ is even, then $(x,y)$ has two neighbours $(x+1,y)$ and $(x,y+1)$ to the right of $(x,y)$ in $\preceq$, and one neighbour $(x-1,y)$ to the left. On the other hand, if $x+y$ is odd, then $(x,y)$ has two neighbours $(x-1,y)$ and $(x,y-1)$ to the left of $(x,y)$ in $\preceq$, and one neighbour $(x+1,y)$ to the right. Consider a vertex $p=(x,y)$ incident to distinct edges $pq$ and $pr$. 
If $p\prec q\prec r$, then $x+y$ is even and $q=(x,y+1)$ and $r=(x+1,y)$, implying that $pq$ and $pr$ and in distinct queues. 
If $r\prec q\prec p$, then $x+y$ is odd and $r=(x-1,y)$ and $q=(x,y-1)$, implying that $pq$ and $pr$ and in distinct queues. 
\end{proof}


\begin{lem}
\label{WallTrack}
The wall has a 4-track layout, such that for all distinct edges $pq$ and $pr$, the vertices $q$ and $r$ are in distinct tracks.
\end{lem}

\begin{proof}
Consider the following vertex ordering of $W$. For vertices $(x,y)$ and $(x',y')$ of $W$, define $(x,y)\preceq (x',y')$ if $x<x'$, or $x=x'$ and $y\leq y'$. 

Colour each vertex $(x,y)$ of $W$ by $(x+2y)\bmod{4}$, as illustrated in \cref{WallFourColouring}. Observe that this defines a proper vertex colouring of $W$. Order each colour class by $\preceq$. Each colour class is now a track. Observe that for all distinct edges $pq$ and $pr$, the vertices $q$ and $r$ are in distinct tracks. Put another way, this is a 4-colouring of the square of $W$. 

\begin{figure}[h]
\begin{center}
\includegraphics{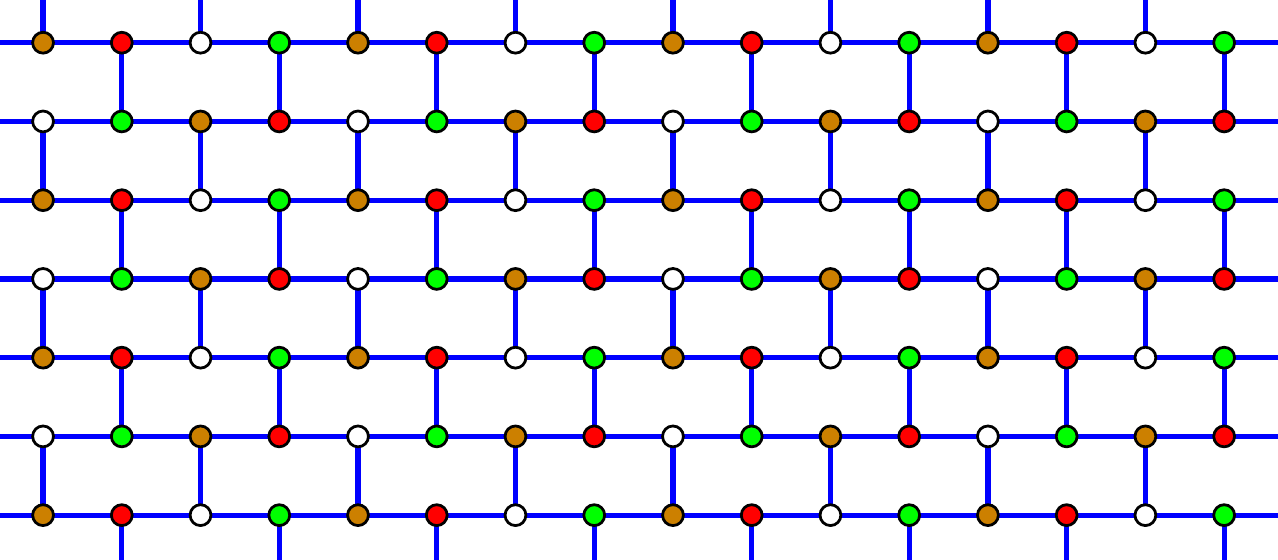}
\caption{4-Colouring the wall}
\label{WallFourColouring}
\end{center}
\end{figure}

Suppose on the contrary that edges $(x_1,y_1)(x_4,y_4)$ and $(x_2,y_2)(x_3,y_3)$ cross, where $(x_1,y_1)\prec(x_2,y_2)$ in some track, and $(x_3,y_3)\prec(x_4,y_4)$ in some other track. Thus $x_1\leq x_2$ and $x_3\leq x_4$. Without loss of generality, $(x_1,y_1)\prec(x_3,y_3)$. Thus 
$(x_1,y_1)\prec(x_3,y_3)\prec(x_4,y_4)$. Hence $x_1\leq x_3\leq x_4$. 

Suppose that $x_1=x_4$. Thus $y_1<y_3<y_4$, implying $y_4\geq y_1+2$ and $(x_1,y_1)(x_4,y_4)$ is not an edge. Now assume that $x_1<x_4$. Since $(x_1,y_1)(x_4,y_4)$ is an edge, $x_4=x_1+1$ and $y_1=y_4$. 

In what follows, all congruences are modulo 4. We have $x_3+2y_3\equiv x_4+2y_4$. Thus $x_3-x_4\equiv 2(y_4-y_3)$, implying $x_3-x_4$ is even. Since $x_1\leq x_3\leq x_4=x_1+1$, we have $x_3=x_4$. 
Since $(x_3,y_3)\prec(x_4,y_4)$, we have $y_3<y_4$. 
Since $x_1+2y_1\equiv x_2+2y_2$, we have $x_1-x_2\equiv 2(y_2-y_1)$, implying $x_1-x_2$ is even. 

Suppose that $x_2\leq x_4$. Then $x_1\leq x_2\leq x_4=x_1+1$. Since $x_1-x_2$ is even, $x_1=x_2$. Since 
$(x_1,y_1)\prec(x_2,y_2)$, we have $y_1<y_2$. Since $(x_2,y_2)(x_3,y_3)$ is an edge and $x_3=x_4=x_1+1=x_2+1$, we have $y_2=y_3$. Similarly, 
since $(x_1,y_1)(x_4,y_4)$ is an edge and $x_4=x_1+1$, we have $y_1=y_4$. Since $y_3<y_4$, we have $y_2<y_1$, which is a contradiction. Now assume that $x_2>x_4$. 

Since $x_2>x_4=x_3$ and $(x_2,y_2)(x_3,y_3)$ is an edge,  $y_2=y_3$. Thus $x_2=x_3+1=x_4+1=x_1+2$. Since $x_1+2y_1\equiv x_2+2y_2$ we have $2y_1\equiv 2 + 2y_2$, implying $y_1-y_2$ is odd. Since $y_4=y_1$ and $y_3=y_2$, we have $y_4-y_3$ is odd. However, since $x_3=x_4$ and $x_3+2y_3\equiv x_4+2y_4$, we have $y_3-y_4$ is even.  This contradiction proves that no two edges between the same pair of tracks cross. 
\end{proof}

\section{The Main Proofs}

Here we give a unified proof of \cref{3Monotone}, \cref{2Queue} and \cref{4Track}. Let $G$ be a two-sided $2k$-monotone bipartite $\epsilon$-expander with bipartition $A,B$. An infinite family of such graphs exist by \cref{TwoSidedBourgain} for fixed $\epsilon$ and $k$. We may assume that $k\geq\frac{3}{\epsilon}$. Let $E_1,\dots,E_{2k}$ be the corresponding partition of $E(G)$. Now define a graph $G'$. For each vertex $v\in A$, introduce the following cycle in $G'$:
\begin{align*}
C_v:=(&\nu_{0},v_{1},\nu_{1},v_{2},\nu_{2},\dots,v_k,\nu_k,\nu_k',v_k',\nu_{k-1}',v_{k-1}',\dots,\nu_1',v_1',\nu_0').
\end{align*}
For each vertex $w\in B$, introduce the following cycle in $G'$:
\begin{align*}
C_w:=(&w_{-1},\omega_{-1},w_0,\omega_{0},w_{1},\omega_{1},w_{2},\omega_{2},\dots,w_k,\omega_k,w_{k+1},\omega_{k+1},w_{k+2},\\
&w_{k+2}',\omega_{k+1}',w_{k+1}',\omega_k',w_k',\omega_{k-1}',w_{k-1}',\dots,\omega_1',w_1',\omega_0',w_0',\omega_{-1}',w_{-1}').
\end{align*}
All the above cycles are pairwise disjoint in $G'$. Finally, for each edge $vw$ of $G$, if $vw\in E_i$ then add the edges $v_iw_i$ and $v_i'w_i'$ to $G'$. 

Observe that $G'$ is bipartite with colour classes:
\begin{align*}
X:=&\{v_i:i\in[1,k]\}\cup\{\nu'_i:i\in[0,k]\}\cup\{w'_i:i\in[-1,k+2]\}\cup\{\omega_i:i\in[-1,k+1]\}\\
Y:=&\{v_i':i\in[1,k]\}\cup\{\nu_i:i\in[0,k]\}\cup\{w_i:i\in[-1,k+2]\}\cup\{\omega'_i:i\in[-1,k+1]\}.
\end{align*}

We now show that \cref{Gen} is applicable to  $G'$. For $v\in A$, let $k_v:=2k+1$. For $w\in A$, let $k_w:=2k+7$. Each cycle $C_v$ has length $2k_v$, as required. 
Since $k\geq\frac{3}{\epsilon}$, we have $\frac{2k+7}{2k+1}\leq 1+\frac{\epsilon}{4}$, as required. 
We now show that the final requirement in \cref{Gen} is satisfied. 
Consider an edge $vw\in E_i$, where $v\in A$ and $w\in B$. Then $v_i\in C_v\cap X$ and $w_i\in C_w\cap Y$ and $v_i'\in C_v\cap Y$ and $w'_i\in C_w\cap X$. Thus  the edges $v_iw_i$ and $v_i'w_i'$ in $G'$ satisfy  the final requirement in \cref{Gen}. 
Hence $G'$ is an $\epsilon'$-expander for some $\epsilon'$ depending on $\epsilon$, $k$ and $\Delta(G)\leq 2k$. 

\begin{figure}[!h]
\begin{center}
\includegraphics[width=\textwidth]{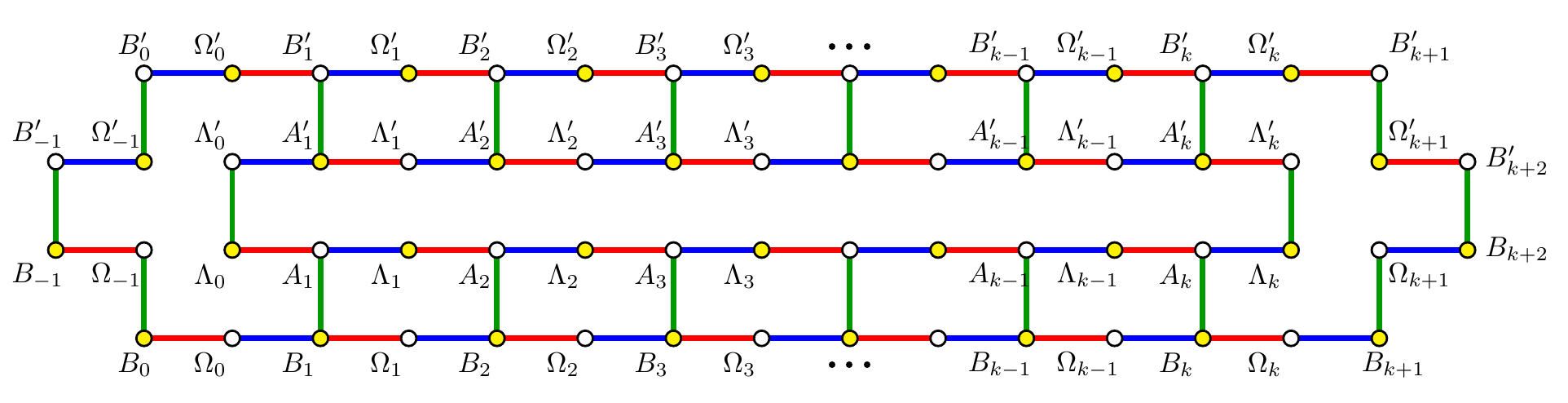}
\caption{The graph $H$. The inner cycle corresponds to vertices in $A$. The outer cycle corresponds to vertices in $B$. \label{SubWall}}
\end{center}
\end{figure}

For $i\in[1,k]$, let $A_i:=\{v_i:v\in A\}$ and $\Lambda_i:=\{\nu_i:v\in A\}$ .
For $i\in[0,k]$, let $A'_i:=\{v'_i:v\in A\}$ and $\Lambda'_i:=\{\nu'_i:v\in A\}$. 
Similarly, for $i\in[-1,k+2]$, let $B_i:=\{w_i:w\in B\}$ and $B'_i:=\{w'_i:w\in B\}$
and $\Omega_i:=\{\omega_i:w\in B\}$ and $\Omega'_i:=\{\omega'_i:w\in B\}$. 
By ordering each of these sets by the given ordering of $A$ or $B$, we consider each such set to be a track. As illustrated in \cref{SubWall}, the graph $H$ obtained from $G'$ by identifying each of these tracks into a single vertex is a subgraph of the  wall.


\begin{figure}[!h]
\begin{center}
\includegraphics[width=\textwidth]{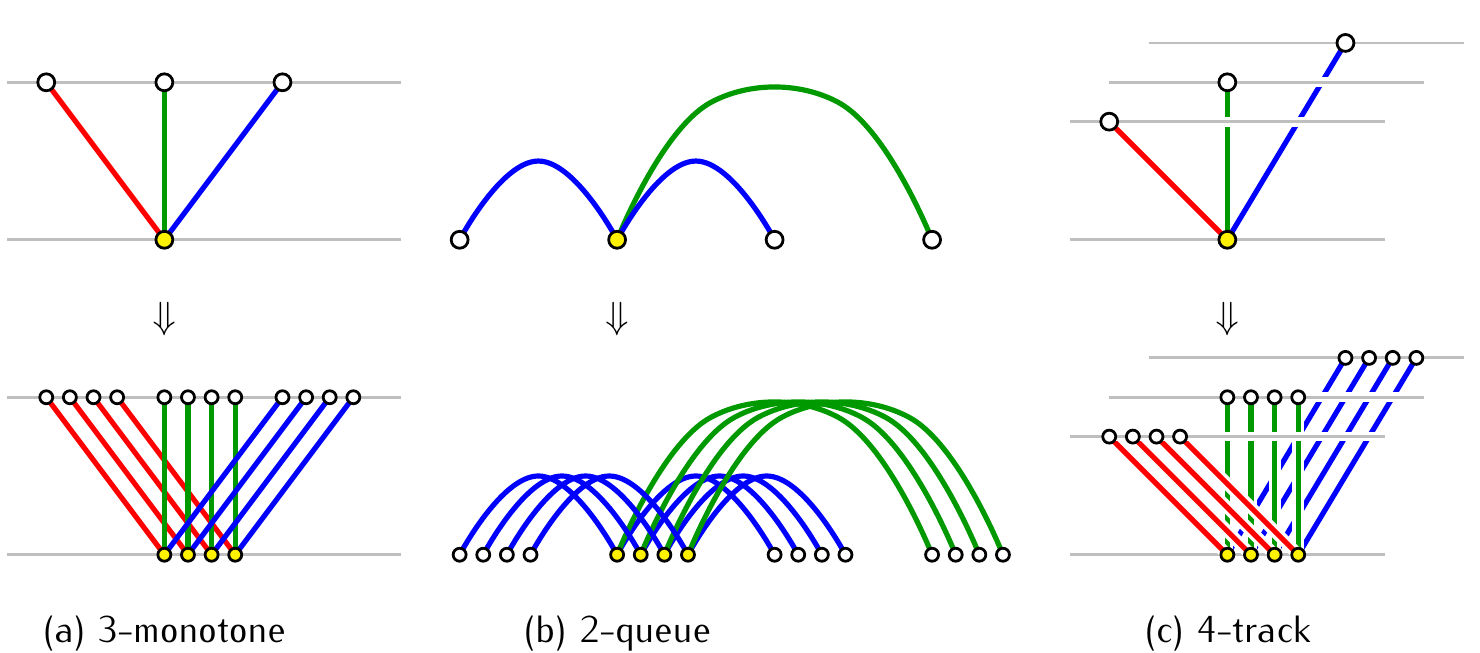}
\caption{Replacing each vertex of $H$ by a track in $G'$.
\label{Uncontract}}
\end{center}
\end{figure}

In other words, there is a homomorphism from $G$ to $H$, where the preimage of each vertex in $H$ is a track in $G$. For each edge $pq$ of $H$, where $pq$ is of the form $a_ib_i$ or $a_i'b_i'$, there is no crossing in $G'$ between the tracks corresponding to $p$ and $q$ since these edges correspond to a monotone matching. For every other edge $pq$ of $H$, the edges between the tracks corresponding to $p$ and $q$ form a non-crossing perfect matching. By \cref{Wall}, $H$ is 3-monotone. Replacing each vertex of $H$ by the corresponding track gives a 3-monotone layout of $G'$, as illustrated in \cref{Uncontract}(a). This proves \cref{3Monotone}. Similarly, by \cref{WallQueue}, $H$ has a 2-queue layout, such that for all edges $pq$ and $pr$ with $p\prec q\prec r$ or $r\prec q\prec p$, the edges $pq$ and $pr$ are in distinct queues. Replacing each vertex of $H$ by the corresponding track gives a 2-queue layout of $G'$, as illustrated in \cref{Uncontract}(b). This proves \cref{2Queue}. Finally, by \cref{WallTrack}, $H$ has a 4-track layout, such that for all edges $pq$ and $pr$, the vertices $q$ and $r$ are in distinct tracks. Replacing each vertex of $H$ by the corresponding track gives a 4-track layout of $G'$, as illustrated in \cref{Uncontract}(c). This proves \cref{4Track}. Note that, in fact, between each pair of tracks, the edges form a monotone matching. 

\section{Open Problems}
\label{OpenProblems}

Heath~\emph{et~al.}~\cite{HR92,HLR92} conjectured that planar graphs have bounded queue-number, which holds if and only if 2-page graphs have bounded queue-number \citep{DujWoo-DMTCS05}. The best upper bound on the queue-number of planar graphs is $O(\log n)$ due to \citet{Duj}; see \citep{DMW13} for recent extensions.

More generally (since planar graphs have bounded page-number \citep{BS84,Yannakakis89}), \citet{DujWoo-DMTCS05} asked whether queue-number is bounded by a function of page-number. This is equivalent to whether 3-page graphs have bounded queue-number \citep{DujWoo-DMTCS05}. 

\citet{DujWoo-DMTCS05}  also asked whether page-number is bounded by a function of queue-number, which holds if and only if 2-queue graphs have bounded page-number \citep{DujWoo-DMTCS05}. 

\citet{GM-JCTB} established a close connection between expanders and linear treewidth that, with \cref{3Page}, gives an infinite family of $n$-vertex $3$-page graphs with $\Omega(n)$ treewidth (and maximum degree 3). This observation seems relevant to a question of \citet{DujWoo11}, who asked whether there is a polynomial time algorithm to determine the book thickness of a graph with bounded treewidth; see \citep{BanEpp} for related results and questions. 

A final thought: 3-page graphs arise in knot theory, where they are called Dynnikov Diagrams \citep{Dynnikov-FAP00,Dynnikov-UMN98,Dynnikov-AAM01}. It would be interesting to see if the existence of 3-page expanders has applications in this domain. 

\subsection*{Note} 

As mentioned earlier, the proof of \cref{Gen} is reminiscent of the replacement product; see \citep{ASS08,RVV-AM02,HLW-BAMS06,DW-ToC10}. A referee observed that it is possible to obtain 3-monotone expanders via the replacement product as follows.  Apply Lemma~3.1 of \citet{DW-ToC10} where $G_1$ is the $2k$-monotone bipartite expander due to \BY, and $G_2$ is a cycle of length $2k$. This gives a 4-monotone bipartite expander (allowing parallel edges). Observe that two of the monotone matchings are the same. Discard one of them to get a 3-monotone graph. The expansion can only drop by a bounded amount, leaving a 3-monotone bipartite expander. 
However, it is unclear whether this construction gives a 2-sided expander as in our construction. Most importantly, the method presented in this paper (using the infinite wall) leads to 2-queue expanders, 4-track expanders, and expanders with simple thickness 2. 

\subsection*{Acknowledgement} 
This research was initiated at the Workshop on Graphs and Geometry held at the Bellairs Research Institute in 2014. Thanks to the referee for the above observation. 


\def\soft#1{\leavevmode\setbox0=\hbox{h}\dimen7=\ht0\advance \dimen7
  by-1ex\relax\if t#1\relax\rlap{\raise.6\dimen7
  \hbox{\kern.3ex\char'47}}#1\relax\else\if T#1\relax
  \rlap{\raise.5\dimen7\hbox{\kern1.3ex\char'47}}#1\relax \else\if
  d#1\relax\rlap{\raise.5\dimen7\hbox{\kern.9ex \char'47}}#1\relax\else\if
  D#1\relax\rlap{\raise.5\dimen7 \hbox{\kern1.4ex\char'47}}#1\relax\else\if
  l#1\relax \rlap{\raise.5\dimen7\hbox{\kern.4ex\char'47}}#1\relax \else\if
  L#1\relax\rlap{\raise.5\dimen7\hbox{\kern.7ex
  \char'47}}#1\relax\else\message{accent \string\soft \space #1 not
  defined!}#1\relax\fi\fi\fi\fi\fi\fi}

\appendix
\section{Separators in Bipartite Expanders}
\label{Separators}

A \emph{separator} in a graph $G$ is a set $Z\subseteq V(G)$ such that each component of $G-Z$ has at most $\frac{|V(G)|}{2}$ vertices. The following connection between expanders and separators is well known, although we are unaware of an explicit proof for bipartite expanders, so we include it for completeness. 

\begin{lem}
\label{ExpanderSeparator}
If $G$ is a bipartite $\epsilon$-expander with $2n$ vertices, then every separator in $G$ has size at least 
$\tfrac{\epsilon}{2}(n-1)-1$. 
\end{lem}

\begin{proof}
Let $A,B$ be the bipartition of $G$ with $|A|=|B|=n$. 
Let $Z$ be a separator of $G$. Our goal is to prove that $|Z|\geq \frac{\epsilon }{2}(n-1)-1$. Let $Z_1:=Z\cap A$ and $Z_2:=Z\cap B$. 

Let $X_1,\dots,X_k$ be a partition of $V(G-Z)$ such that each $X_i$ is the union of some subset of the components of $G-Z$ with at most $n$ vertices in total, and subject to this condition, $k$ is minimal. This is well-defined, since each component of $G-Z$ has at most $n$ vertices. By minimality, $|X_i|+|X_j|>n$ for all distinct $i,j\in[1,k]$. If $k\geq 4$ then $|X_1|+|X_2|>n$ and $|X_3|+|X_4|>n$, which contradicts the fact that $|V(G)|=2n$. Hence $k\leq 3$. Let $A_i:=X_i\cap A$ for $i\in[1,k]$. 

First suppose that $|A_i|\geq\frac{n}{2}$ for some $i\in[1,k]$. Let $S$ be a subset of $A_i$ with exactly $\floor{\frac{n}{2}}$ vertices. Observe that $N(S)\subseteq (X_i\setminus S)\cup Z$. Thus $$(1+\epsilon)|S|\leq |N(S)|\leq|X_i| -|S|+|Z|\leq n-|S|+|Z|,$$
and
$$|Z|\geq (2+\epsilon)|S|-n = (2+\epsilon)\FLOOR{\frac{n}{2}}-n\geq 
(2+\epsilon)\left(\frac{n-1}{2}\right)-n
=\frac{\epsilon}{2}(n-1)-1,$$
as desired. 

Now assume that $|A_i|<\frac{n}{2}$ for all $i\in[1,k]$. Observe that 
\begin{align*}
\sum_i(2+\epsilon)|A_i|-|X_i|
=&\big((2+\epsilon)\sum_i|A_i|\big)-\big(\sum_i|X_i|\big)\\
=&(2+\epsilon)(n-|Z_1|)-(2n-|Z_1|-|Z_2|)\\
=&\epsilon n -(1+\epsilon)|Z_1|+|Z_2|.
\end{align*}
Thus, for some $i\in[1,k]$, we have $(2+\epsilon)|A_i|-|X_i|\geq \frac1k (\epsilon n -(1+\epsilon)|Z_1|+|Z_2|)$. 
Observe that $N(A_i)\subseteq (X_i\setminus A_i)\cup Z_2$. 
Thus $$(1+\epsilon)|A_i|\leq |N(A_i)|\leq|X_i| -|A_i|+|Z_2|,$$
and
$$|Z_2|\geq (2+\epsilon)|A_i|-|X_i|\geq \tfrac1k (\epsilon n -(1+\epsilon)|Z_1|+|Z_2| ), $$
implying
$$(k-1)|Z_2| + (1+\epsilon)|Z_1| \geq \epsilon n. $$
Since $k\leq 3$ and $1+\epsilon\leq 2$, we have $2|Z| \geq (k-1)|Z_2| + (1+\epsilon)|Z_1| \geq \epsilon n$, implying $|Z|\geq\frac{\epsilon n}{2}$ as desired. 
\end{proof}

\section{Subdivisions}
\label{SubdivExpander}

Here we show that the 2-subdivision of a bipartite expander is another bipartite expander. This result is well known, although we are unaware of an explicit proof, so we include it for completeness. 

\begin{lem}
For every two-sided bipartite $\epsilon$-expander $G$ with maximum degree $d$, if  $G'$ is the graph obtained from $G$ by subdividing each edge twice, then $G'$ is a two-sided bipartite $\epsilon'$-expander, for some $\epsilon'$ depending on $\epsilon$ and $d$. 
\end{lem}

\begin{proof}
Say $G$ has $m$ edges, and $(A,B)$ is the bipartition of $G$ with $n=|A|=|B|$. Since $G$ is an $\epsilon$-expander, each vertex has degree greater than 1 (and at most $d$). Thus $2n\leq m\leq dn$. Observe that $G'$ is bipartite with bipartition $(A\cup A',B\cup B')$, where for each edge $e$ of $G$, 
exactly one division vertex of $e$ is in $A'$, and exactly one division vertex of $e$ is in $B'$. Each colour class of $G'$ has $n+m$ vertices. 

Let $S\subseteq A$ and $S'\subseteq A'$ such that $|S|+|S'|\leq \half(n+m)$. By the symmetry between $A'$ and $B'$, it suffices to prove that $|N(S\cup S')|\geq(1+\epsilon')|S\cup S'|$, for some $\epsilon'$ depending solely on $\epsilon$ and $d$. We do so with the following definition of $\epsilon'$:
\begin{align*}
\beta&:=\frac{4d+3}{4d+4}\\
\gamma&:=\frac{\epsilon(1-\beta)}{2(\epsilon+\beta)} \\
\epsilon'&:=\min\left\{(1+\epsilon)(1-\gamma)-1,\,
\frac{(\epsilon+\beta)(1-\gamma)}{(1+\epsilon)\beta}-1,\,
\frac{\gamma}{d}\right\}.
\end{align*}
We now show that $\epsilon'>0$. 
Since $1-\epsilon<1<2\beta$, we have 
$1-\beta<\epsilon+\beta$ and $\frac{1-\beta}{\epsilon+\beta}<1$. Thus $0<\gamma<\frac{\epsilon}{2}<1$. 
Since $\gamma<\frac{\epsilon}{2}$, we have $(1+\epsilon)(1-\gamma)>1$, and since
$\gamma<\frac{\epsilon(1-\beta)}{(\epsilon+\beta)}$, it follows that 
$\frac{(\epsilon+\beta)(1-\gamma)}{(1+\epsilon)\beta}>1$. [\emph{Proof}. 
$\frac{\epsilon(1-\beta)}{(\epsilon+\beta)} >\gamma$
implies
$\epsilon(1-\beta) >\gamma(\epsilon+\beta)$
implies
$(\epsilon+\beta)-(1+\epsilon)\beta >\gamma(\epsilon+\beta)$
implies 
$(\epsilon+\beta)(1-\gamma)>(1+\epsilon)\beta$
implies
$\frac{(\epsilon+\beta)(1-\gamma)}{(1+\epsilon)\beta}>1$.]\ Thus $\epsilon'>0$. Also note that $\epsilon'
\leq\frac{\gamma}{d}\leq\frac{\epsilon}{2d}\leq\frac{1}{2d}\leq 1$.

Consider a subdivided edge $(v,w',v',w)$ of $G$ where $v\in A$, $w\in B$, $v'\in A'$ and $w'\in B'$. 
Say $vw$ is \emph{type-1} if $v\in S$ and $v'\in S'$. Let $m_1$ be the number of type-1 edges in $G$. 
Say $vw$ is \emph{type-2} if $v\in S$ and $v'\not\in S'$. Let $m_2$ be the number of type-2 edges in $G$. 
Say $vw$ is \emph{type-3} if $v\not\in S$ and $v'\in S'$. Let $m_3$ be the number of type-3 edges in $G$. 

Let $X$ be the set of vertices in $B$ adjacent in $G$ to some vertex in $S$ and adjacent in $G'$ to some vertex in $S'$. 
The endpoint of each type-1 edge in $B$ is in $X$, and each vertex in $X$ has degree at most $d$. Thus $m_1\leq d|X|$. 
Let $Y$ be the set of vertices in $B$ adjacent in $G$ to some vertex in $S$ and adjacent in $G'$ to no vertex in $S'$.
Each vertex in $Y$ is incident to some type-2 edge, and distinct vertices in $Y$ are incident to distinct type-2 edges. Thus $m_2\geq |Y|$. 
Let $Z$ be the set of vertices in $B$ adjacent in $G$ to no vertex in $S$ and adjacent in $G'$ to some vertex in $S'$.
The endpoint of each type-3 edge in $B$ is in $X\cup Z$, and each vertex in $X\cup Z$ has degree at most $d$. Thus $m_1+m_3\leq d(|X|+|Z|)$. 

Note that $X,Y,Z$ are pairwise disjoint. By the definition of $X$ and $Y$, we have $N_G(S)=X\cup Y$. 
Observe that $|N_{G'}(S\cup S')\cap B'|=m_1+m_2+m_3$ and $N_{G'}(S\cup S')\cap B=X\cup Z$. Thus 
\begin{equation}
\label{NG'SS'}
|N_{G'}(S\cup S')|=|X|+|Z|+m_1+m_2+m_3.
\end{equation}
 Also note that $|S'|=m_1+m_3$. Since $|S|+|S'|\leq\half(n+m)$, 
\begin{equation}
\label{m1m3}
m_1+m_3 \leq \half(n+m) - |S|.
\end{equation}



First suppose that $m_2\geq(1+\epsilon')\frac{m}{2}+m_3$. 
Thus $m_2\geq(1+\epsilon')n+m_3\geq (1+\epsilon')|S|+m_3$. 
Since $|X|\geq\frac{m_1}{d}\geq\epsilon' m_1$ and $1\geq\epsilon'$, we have
$m_2+|X|\geq(1+\epsilon')|S|+\epsilon'(m_1+m_3)$. 
By \eqref{NG'SS'}, 
$$|N(S\cup S')|=m_1+m_2+m_3+|X|+|Z|\geq(1+\epsilon')(|S|+m_1+m_3)=(1+\epsilon')|S\cup S'|,$$
as desired. 

%
Now assume  that $m_2\leq(1+\epsilon')\frac{m}{2}+m_3$. Since $m\leq dn$ and $\epsilon'\leq\frac{1}{2d}$, 
\begin{equation}
\label{blah}
2m_2-m-2m_3\leq\epsilon'm \leq\tfrac{n}{2}.
\end{equation}
Since $G$ contains $m-m_1-m_2$ edges incident to $A-S$, and each vertex in $A-S$ has degree at most $d$, we have $m-m_1-m_2\leq d(n-|S|)$. By \eqref{m1m3}, 
$$m-m_2\leq d(n-|S|)+m_1\leq  d(n-|S|)+\half(n+m)-|S|-m_3.$$
By \eqref{blah},
\begin{equation}
\label{blahblah}
(2d+2)|S|\leq   (2d+1)n+2m_2-m-2m_3\leq   (2d+1)n+\tfrac{n}{2}.
\end{equation}
Therefore $|S|\leq\beta n$. Say $|S|=\alpha n$, where $0<\alpha\leq \beta$.

If $|S|\leq\frac{n}{2}$ then (since $G$ is an $\epsilon$-expander), $$|X|+|Y|=|N_G(S)|\geq(1+\epsilon)|S|\geq\frac{1+\epsilon'}{1-\gamma}\,|S|.$$ 
We can reach the same conclusion when  $|S|\geq\frac{n}{2}$ as follows. Considering any subset of $S$ of size $\floor{\frac{n}{2}}$, we have $|X|+|Y|=|N_G(S)|\geq(1+\epsilon)\floor{\frac{n}{2}}>\frac{n}{2}$. Thus $|B-(X\cup Y)|\leq\frac{n}{2}$. Since $G$ is a two-sided $\epsilon$-expander,   
$|N_G(B-(X\cup Y))|\geq(1+\epsilon)|B-(X\cup Y)|=(1+\epsilon)(n-|X|-|Y|)$. No vertex in $B-(X\cup Y)$ is adjacent to $S$. Thus
$$n-|S|\geq |N_G(B-(X\cup Y))|\geq(1+\epsilon)(n-|X|-|Y|).$$ That is,
\begin{equation}
\label{XY}
|X|+|Y|\geq \frac{\epsilon n+|S|}{1+\epsilon} = \frac{\epsilon+\alpha}{1+\epsilon}\,n.
\end{equation}
Since $\alpha \leq \beta$, it follows that $\frac{\epsilon+\beta}{\beta} \leq \frac{\epsilon+\alpha}{\alpha}$. 
Thus
$$1+\epsilon' \leq \frac{(\epsilon+\beta)(1-\gamma)}{(1+\epsilon)\beta} \leq \frac{(\epsilon+\alpha)(1-\gamma)}{(1+\epsilon)\alpha}.$$
Hence
$$\frac{(1+\epsilon')\alpha}{1-\gamma} \leq \frac{\epsilon+\alpha}{1+\epsilon}.$$
By \eqref{XY}, 
$$|X|+|Y|\geq \frac{\epsilon+\alpha}{1+\epsilon}\,n\geq \frac{(1+\epsilon')\alpha}{1-\gamma}\,n=\frac{1+\epsilon'}{1-\gamma}\,|S|,$$
as claimed. 

Since $m_2\geq |Y|$ and $0<\gamma<1$, 
$$(1+\epsilon')|S|\leq(1-\gamma)( |X|+|Y|)\leq(1-\gamma) |X|+m_2\leq(1-\gamma)( |X|+|Z|)+m_2.$$
Since $\epsilon'\leq\tfrac{\gamma}{d}$ and $m_1+m_3\leq d(|X|+|Z|)$, we have $\epsilon'(m_1+m_3)\leq\gamma(|X|+|Z|)$. Hence
$$(1+\epsilon')|S|+\epsilon'(m_1+m_3)\leq |X|+|Z|+m_2.$$
Therefore, by \eqref{NG'SS'}, 
$$(1+\epsilon')|S\cup S'|=(1+\epsilon')(|S|+m_1+m_3)\leq |X|+|Z|+m_1+m_2+m_3=|N_{G'}(S\cup S')|.$$
This completes the proof. 
\end{proof}

\end{document}